\documentclass{amsart}
\usepackage{graphicx}
\usepackage{amsfonts}

\newcommand{\R}{{\Bbb R}}
\newtheorem{theorem}{Theorem}[section]

 \newtheorem{lemma}[theorem]{Lemma}
\newtheorem{corollary}[theorem]{Corollary}
\newtheorem{conjecture}[theorem]{Conjecture}

\theoremstyle{definition}

\newtheorem{example}[theorem]{Example}

\theoremstyle{remark}
\newtheorem{remark}[theorem]{Remark}

\numberwithin{equation}{section}
\begin{document}

\title[$3/2$ stability condition]{Wright type delay  
differential equations with negative Schwarzian}

\thanks{This research was supported by FONDECYT (Chile), project 8990013.
E. Liz was supported in part by D.G.E.S. (Spain), project
PB97-0552. V. Tkachenko was supported in part by F.F.D. of Ukraine,
project 01.07/00109.}

\author[E. Liz, M. Pinto, G. Robledo, V. Tkachenko, S. Trofimchuk]{}
\email{eliz@dma.uvigo.es}
\email{pintoj@abello.dic.uchile.cl}
\email{gonzalo.robledo@etud.univ-pau.fr}
\email{trofimch@uchile.cl}
\email{vitk@imath.kiev.ua}


\date{\today}

\maketitle

\centerline{\scshape  Eduardo Liz}
\medskip

{\footnotesize
\centerline{ Departamento de Matem\'atica Aplicada, E.T.S.I. Telecomunicaci\'on }
\centerline{ Universidad de Vigo, Campus Marcosende, 36280 Vigo, Spain}
}
\medskip
\centerline{\scshape  Manuel Pinto, Gonzalo Robledo, and Sergei Trofimchuk}
\medskip

{\footnotesize
\centerline{ Departamento de Matem\'aticas, Facultad de Ciencias }
\centerline{ Universidad de Chile, Casilla 653, Santiago, Chile}
}

\medskip
\centerline{\scshape  Victor Tkachenko}
\medskip

{\footnotesize
\centerline{ Institute of Mathematics,
National Academy of Sciences of Ukraine}
\centerline{ Tereshchenkivs'ka str. 3, Kiev, Ukraine}
}

\bigskip
\begin{quote}{\normalfont\fontsize{8}{10}\selectfont
{\bfseries Abstract.}
 We prove that the well-known  
$3/2$ stability condition established for the Wright equation (WE) 
still holds if
the nonlinearity $p(\exp(-x)-1)$ in WE is replaced by a decreasing or 
unimodal smooth function $f$ with $f'(0)<0$ satisfying
the standard negative feedback and below boundedness conditions
 and having everywhere  negative 
Schwarz derivative. \par}
\end{quote}

\section{Introduction}
\label{sec.int}
In this paper we study the global stability properties 
of the scalar delay-differential equation
\begin{equation}
x'(t) = f(x(t-1)), \label{wr}
\end{equation}
where $f \in C^{3}(\R,\R)$ satisfies the following 
additional conditions \textbf{(H)}: 

{\bf (H1)} $xf(x) < 0$  for  $x\neq 0$ and $f'(0) < 0$.

{\bf (H2)} $f$ is bounded below 
and has at most one critical point $x^*\in \R$ which is a local extremum.

{\bf (H3)} $(Sf)(x)<0$ for all $x\neq x^*$, where
$
Sf=f'''(f')^{-1}-3/2(f'')^2(f')^{-2}
$
is  the  Schwarz derivative of $f$. 

The negative feedback condition {\bf (H1)} and boundedness condition 
{\bf (H2)} are very typical in the theory of (\ref{wr}); the first one
causes solutions to tend to oscillate about zero, while both of them 
guarantee the existence of the global compact attractor to 
Eq. (\ref{wr}) (e.g. see \cite{mp}). On the other hand, 
the Schwarzian negativity condition {\bf (H3)} is rather common 
in the theory of one-dimensional dynamical systems 
(see \cite{stri, sh}) while it is not very usual for the 
studies of delay-differential equations. Here, we introduce {\bf (H3)} 
in  hope of obtaining an analogue of the Singer global 
stability result for an one-dimensional map $f:I\to I$; this result  (see its complete
formulation below)  states that
the local stability of a unique fixed point $e \in I$ plus an appropriate 
(monotone or unimodal) form of the map $f$ imply the 
global attractivity of the equilibrium $e \in I$ \cite{stri}. 
The negative Schwarzian condition is not artificial at all, it  appears naturally also in some
other contexts of the theory of delay differential equations, see e.g.
\cite[Sections 6--9]{mn}, where it is explicitely used, and  \cite[Theorem 7.2, p.~388]{hl},
where the condition
$Sf < 0$ is implicitely required; moreover,
many nonlinear delay differential equations used 
in mathematical modelling 
in biology (e.g. Mackey-Glass, Lasota-Wazewska, Nicholson, 
Goodwin equations) have their right-hand sides satisfying 
the hypothesis {\bf (H3)}. Take also, for example, the celebrated
Wright equation which was used to describe the distribution of prime numbers or 
to model population dynamics of a single species: 
\begin{equation}
\label{we}
x'(t) = - px(t-1)[1+x(t)], \quad p >0.
\end{equation}
For $x(t)>-1$, Eq. (\ref{we}) is reduced to (\ref{wr})
after applying the transformation $y = -\ln(1+x)$. 
In this case $f(x) = (\exp(-x) - 1)$ and, by abuse of  notation, 
we will again  refer to the transformed system
\begin{equation}
x'(t) = p(\exp(-x(t-1)) - 1), \ p > 0,
\label{wri}
\end{equation}
as  the Wright equation. In this case, $f$ is strictly decreasing
and has no  inflexion points; both these facts simplify considerably 
the investigation of (\ref{wri}). Below, we will present two other  
important examples, with nonlinearities which may have an inflexion point
(some ``food-limitation" models) or even a local extremum (population model
exhibiting the Allee effect). 

Eq. (\ref{wr}) has been considered before by several authors
but only assuming conditions {\bf (H1)} and {\bf (H2)},  see \cite{mp, mps, walt, wal} and references therein.
In particular, the Morse decomposition of its 
compact attractor has been described in detail in \cite{mp}. Moreover, 
it was proved also that the Poincar\'e-Bendixson 
theorem holds for (\ref{wr}) with the decreasing nonlinearity $f$
so that the asymptotic periodicity 
is the ``most complicated" type of behavior 
in (\ref{wri}) \cite{mp, mps, wal}. It should be noted also that 
the above information has essentially a ``qualitative" character.  
So that adding {\bf (H3)}, we can hope to obtain some additional 
information of analytic nature about possible 
bifurcations in parametrized families of (\ref{wr}).
The following variational equation along 
its unique steady state $x=0$ plays a very important role in the 
study  of such bifurcations: 
\begin{equation}
x'(t) = f'(0)x(t-1). 
\label{ve}
\end{equation}
As  is well-known, this equation is unstable 
when $-f'(0) > \pi/2$, and  this instability implies the 
existence of  slowly oscillating periodic solutions to (\ref{wr}) 
(see e.g. \cite{walt}). Surprisingly,  the dynamically  simpler case 
$-f'(0) < \pi/2 = 1.571...$ has not been studied thoroughly before, and, 
in particular, it seems that the following
Wright conjecture has not been solved: the inequality $p < \pi/2$ is sufficient
for the global stability in (\ref{wri}). On the other hand, 
the sufficiency of the stronger condition $p < 3/2$ for the 
global stability of Eq. (\ref{wri}) was proved 
in `a very difficult theorem of Wright \cite{WR}'
(see \cite[page 64]{HMO}), where also the  sharper conditions 
$p < 37/24 = 1.5416...$ and $p < 1.567...$ were announced. 
It should be noted that proofs of the $3/2$ stability
condition for Eq. (\ref{wri}) have strongly used  the specific 
exponential form of the nonlinearity $f(x) = p(\exp(-x)-1)$ and, 
in particular, the monotonicity properties of such $f$. 
This fact explains why any analogue
of this Wright result has  not been proved for other,  
essentially different (nonexponential), right-hand 
sides in Eq. (\ref{wr}) (even for monotone $f$, the general
situation being considerably more complicated). 

An important step in solving the Wright 
conjecture was made in Theorem 3 from \cite{walt} which
provides us with some examples of Eq. (\ref{wr}) which satisfy
{\bf (H1), (H2)} and have slowly oscillating periodic solutions, 
even when the corresponding linearized equation (\ref{ve}) 
is exponentially stable. This {means} that 
the local exponential stability 
of the steady state in (\ref{wr}) (or, what is the same, 
exponential stability in (\ref{ve})) with $f$ having only
these standard and usual properties {\bf (H1), (H2)} in general 
does not imply the global asymptotic stability in (\ref{wr}).  
Moreover, as a simple consequence of an elegant approach towards 
stable periodic orbits for scalar equations of the form 
$x'(t) = - \mu x(t) + f(x(t-1)),\ \mu \geq 0$ with Lipschitz 
nonlinearities proposed recently in \cite{waln} (see also \cite{waltdcds}), we get the 
following 
\begin{theorem}[\cite{waln}]
\label{11}
For every $\alpha \geq 0$ 
there exists a smooth strictly decreasing function $f(x)$ satisfying 
{\bf (H1)}, {\bf (H2)}, $-f'(0) = \alpha$ and such that Eq. 
(\ref{wr}) has a nontrivial periodic solution which is hyperbolic,
stable and exponentially attracting with asymptotic phase
(so therefore (\ref{wr}) is  not globally stable). 
\end{theorem}
This result is of  special importance for us, since it shows clearly that 
the strong correlation between local (at zero) and global properties 
of Eq. (\ref{wri}) can not 
be explained only with the concepts presented in 
{\bf (H1), (H2)}.

On the other hand,  Walther's result from \cite{walt} 
cannot be  applied to Eq. (\ref{wri}) so that the original
Wright conjecture remains open. 
We explain here this particularity of Eq. (\ref{wri})
 by its additional
property of having negative 
Schwarz derivative $Sf$: in fact, no function from  \cite[Theorem 3]{walt} 
can have $Sf < 0$. Moreover, bearing in mind   the  following result of D. Singer  
for one-dimensional maps: \textit{
``Assume that the function $h\in C^3[a,b]$
is either strictly decreasing or  has only one critical point $x^*$ (local extremum) 
in $[a,b]$. If $h$ has a unique fixed point
$e\in [a,b]$ which  is locally asymptotically stable and
$(Sh)(x)<0$ for all $x\neq x^*$, then $e$ is the global attractor of the dynamical system $h:[a,b]\to [a,b]$"},
we propose to generalize  Wright's conjecture in the form stated below.
\begin{conjecture}
\label{conj}
 Let all conditions {\bf (H)} be satisfied and 
$-f'(0) < \pi/2$. Then $\lim\limits_{t\to +\infty}x(t)=0$ for every solution $x(t)$ to Eq. (\ref{wr}).
\end{conjecture}
We remark that this conjecture
is very close to the Hal Smith conjecture \cite{smith} to the effect
that the  local and global asymptotic stabilities 
for Nicholson's blowflies equation
$$x'(t) = -\delta x(t) + p x(t-1)\exp(-a(x(t-1))), \quad x,\delta, p, a  > 0, $$
are equivalent.
Observe that this equation  also has a unique positive steady state and  nonlinearity 
satisfying the negative feedback and the negative Schwarzian conditions 
(see \cite{gyt, delta, sya} for  further discussions).

Furthermore, due to recent results of Krisztin  \cite{Kr}, now we
can indicate some class of symmetric and monotone nonlinearities
(e.g. $f(x) = - p\tanh x, f(x) = -p \arctan x, \ p > 0$) for which
the above conjecture is true. Although for  both the  mentioned
functions condition $Sf < 0$ holds, in general the additional convexity condition
imposed  on $f$ in \cite{Kr} (see also \cite{KW}) is
different from {\bf (H3)}: evidently, the requirement of the
negative Schwarzian is not the unique way to approach the problem
(the same situation that we have in the theory of one-dimensional
maps). In fact, to prove our main result (Theorem \ref{5678}), we
only need  some geometric consequences of the inequality $Sf < 0$
for the graph of $f$. For instance, if $f''(0)=0$ then this
geometric consequence is given by $(f(x) - f'(0)x)x > 0$ for
$x\not=0$ (that necessarily holds also under the above mentioned convexity,
symmetry and monotonicity assumptions from \cite{Kr,KW}). This
geometric approach was developed further in \cite{delta}, where a
generalization of the Yorke condition \cite[Section 4.5]{kuang} was proposed
instead of {\bf (H3)}.

In this paper we carry out the first step towards  Conjecture \ref{conj}, 
showing how all the conditions ({\bf H}) come together  in proving 
\begin{theorem}
\label{5678}
If, in addition to  {\bf (H)}, we have  $-f'(0) \leq 1.5$, 
then the steady state 
solution $x(t) = 0$ of Eq. (\ref{wr}) is globally
attracting. 
\end{theorem}
To prove Theorem \ref{5678}, we will essentially use  an idea
from \cite{ilt}, which allows us to construct some 
one-dimensional map inheriting some
properties of Eq. (\ref{wr}). Roughly speaking, we
consider maps $F_k = F_k(z): \R \to \R, \ F_k(0) = 0$, which give the
values of the $k$-th consecutive extremum of the oscillating solutions
$x(t, z), \ z \not = 0,$ satisfying $x(s,z) \equiv z, \ s \in
[-1,0]$.  Then we investigate some relations existing between the
global attractivity properties of $F_k$ and (\ref{wr}),  trying to deduce in 
this way the
global asymptotical stability of (\ref{wr}) from the corresponding property of
the discrete dynamical system generated by $F_k$.  Since the computation
becomes more and more complicated with the growth of $k$, we only consider the
simplest case $k =1$ here. Computer experiments show that, increasing $k$,
we obtain better approximations to the condition $-f'(0) \leq
\pi/2$ given in Conjecture \ref{conj} (for example, for $k=2$ we
get $-f'(0) \leq 37/24$ and so on). However, due to the technical
complications, this way to approach the above conjecture could be
used only for very special cases.

Curiously, the
note \cite{ilt} devoted to the study of the Yorke type 
functional differential equations with sublinear nonlinearity 
(see \cite{yorke}) and, in particular, the Yorke 3/2 stability criterion, 
still can be extended to the class of nonlinear Wright's type 
delay differential equations. 
We consider the special nature of the number $3/2$ 
(which was found almost simultaneously by A.D. Myshkis \cite{my} 
and by E.M. Wright \cite{WR}) as an
invariant of such a prolongation. Moreover,  
there exists an interesting interplay between both these 
types of functional differential equations if we  consider
a variable coefficient $p(t)$ instead of the constant $p$ in 
Eq. (\ref{wri}) (see \cite{sy} for details). 
In any case, it should be noted that the Yorke and the 
Wright type equations have rather different  structures
(see \cite{ilt,sy,it} for more comments). 

Completing our discussion, we consider briefly
two other Wright type equations studied recently by several authors:

\begin{example}The ``food-limitation" model \cite[p.~456]{sy}
or, what is basically the same, the Michaelis-Menten single species growth equation
with one delay (see  \cite[p.~132]{kuang}):
\begin{equation}
\label{mm}
x'(t) = -r(1+x(t))\frac{x(t-h)}{1+cr(1+x(t-h))}, \ x > -1,\ c\geq 0,\,r, h >0\ . 
\end{equation}
Note that (\ref{mm}) is of the form 
$x'(t) = -r(1+x(t))g(x(t-h))$ with $Sg=0$. 
The change of variables $x = \exp(-y)-1$ reduces (\ref{mm}) 
to $y'(t) = rf(y(t-h))$, where $f(y)=g(\exp(- y)-1)$ is strictly 
decreasing with $f(0) =0,\ f'(0)= -(1+cr)^{-1},$ and $(Sf)(y)= -1/2 <
0$ for all $y$.  By Theorem \ref{5678}, 
the inequality $rh\leq (3/2)(1+ cr)$ implies 
the global stability of the zero solution to Eq. (\ref{mm})
(compare with \cite{kuang} and \cite{sy}). We also point out that Eq. (\ref{mm})
with $c=0$ coincides with the Wright equation (\ref{we}), so that Wright's $3/2$ stability
theorem is a very special case of our result.
\end{example}

\begin{example} Consider now a population 
model described by 
\begin{equation}
\label{ale}
x'(t) = x(t)[a+ bx(t-h) - cx^2(t-h)], \ a, c \in (0,+\infty),\ b \in \R.
\end{equation}
This equation has a unique positive equilibrium 
$e^*$ and the change of variables $x = e^*\exp(-y)$ transforms 
(\ref{ale}) into (\ref{wr}) with $f(x) = -(a+be^*\exp(-x) - c(e^*\exp(-x))^2)$.
Note that $f$ has negative S-derivative as a composition of a quadratic 
polynomial and the real exponential function. If $b \leq 0$ then $f$ is strictly 
decreasing, and if $b>0$, then $f$ has exactly one critical point (minimum).
In the latter case, the population model exhibits the so-called 
Allee effect \cite{mur}. Applying Theorem \ref{5678}, 
we see that $e^*$ attracts all positive solutions of (\ref{ale}) once 
$(2ce^*-b)he^* \leq 1.5$ (compare e.g.  with \cite[pp.~143-146]{kuang},
where also other references can be found). 
\end{example}
The paper is organized as follows.  In Section 2 we define several 
auxiliary scalar functions and study their properties as well as 
relations connecting them. Finally, in the last section, we 
use these functions to prove Theorem \ref{5678}
(notice that, in contrast with \cite{gyt}, we may 
not   use  the above formulated Singer's  result for this purpose). 

\section{Auxiliary functions}

To prove our main result, we will proceed in analogy to \cite{ilt}, 
so that the construction of one-dimensional maps inheriting 
attractivity properties of the dynamical system generated by 
Eq. (\ref{wr}) is the main tool here. In this section, we introduce 
several such scalar maps and study their properties as well 
as the relations existing among them. 

First we note that we can only have eight different possibilities 
for the maps satisfying hypotheses {\bf (H)},  depending on the situation of
the eventual critical point and the inflexion points. A graphic representation of all these
cases is given in Fig. 1 below, where $x^*$ denotes the critical point and $c_1,\ c_2$ are the
inflexion points. We recall that a real function has at most one inflexion point in any
interval in which the  Schwarz derivative is well defined and is negative (see   \cite{sh}).
\begin{figure}[ht]
\centering
\includegraphics{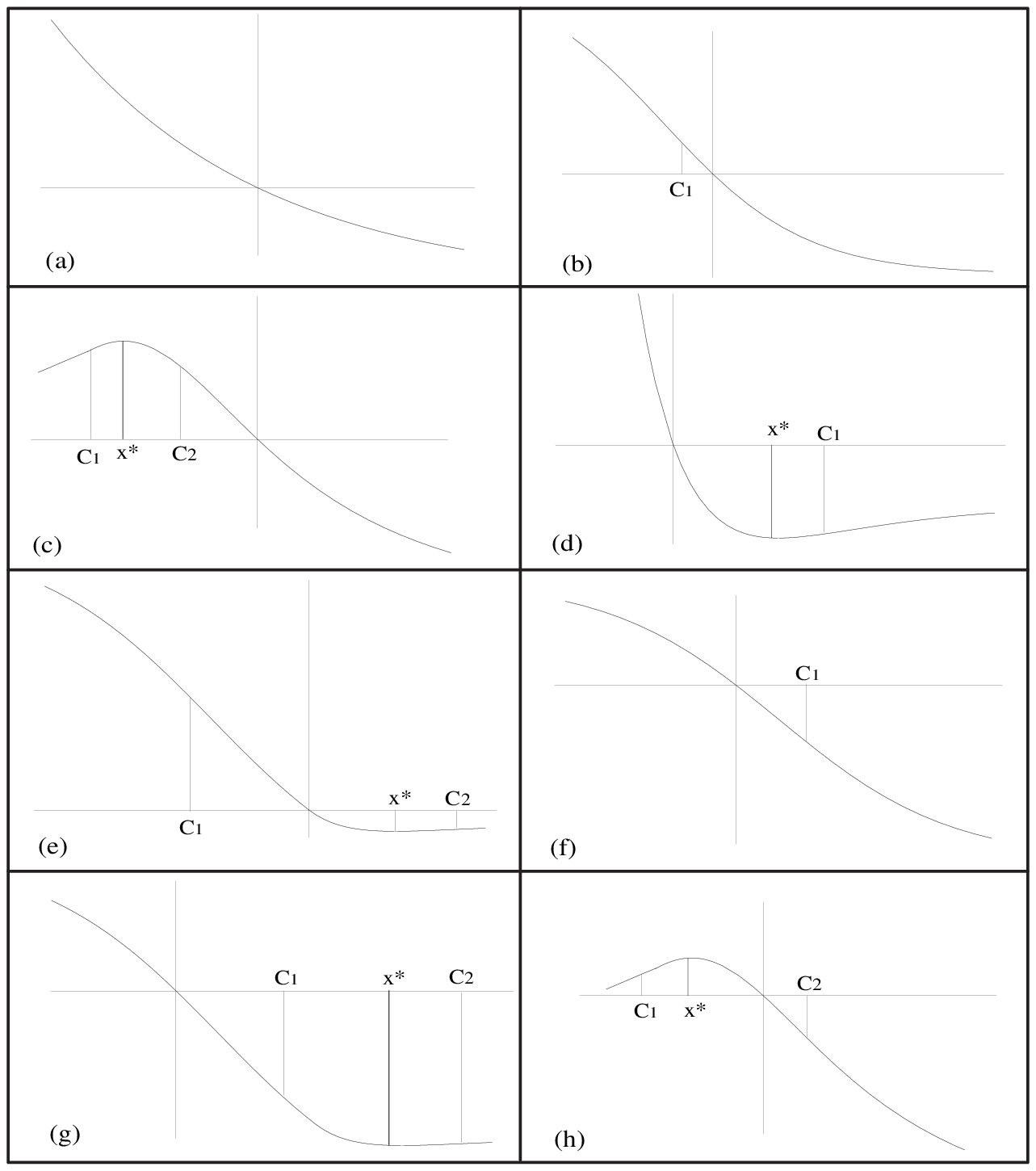}
\caption{Different maps satisfying {\bf (H)}} \label{fig}
\end{figure}

In the sequel, up to the proof of Theorem 
\ref{5678} and with the unique exception made for Corollary \ref{bou}, 
we will always assume that $f$ satisfies {\bf (H)}  and $f''(0) > 0$.  
This situation corresponds to the pictures (a)-(e) 
in Fig. 1.

Next, for $a<0$, $b>0$, we introduce
the set 
\begin{eqnarray*}
K_{a,b}  &= & \{y \in C^{3}(\R):  y\ \mbox{\rm satisfies {\bf (H1),\ (H2)}}, \
 y'(0)= a,\ y''(0) = 2b, \\ & & (Sy)(x) \leq  0   
{\rm\ for \ all \ } x \ {\rm such\ that\ } y'(x) \not= 0  \}.
\end{eqnarray*}
Also, for every $a<0$ and every $b >0$, consider the rational function 
$r(x,a,b)= a^2x/(a-bx)$ defined over $(ab^{-1},+\infty)$.
Let $K_{a,b}^+$ (respectively, $K_{a,b}^-$) be the set 
of  restrictions of  elements of $K_{a,b}$ to   $[0,+\infty)$ (respectively,
to   $(ab^{-1},0]$).   We will denote by $r^+$ and $r^-$ the  restrictions of
$r(\cdot ,a,b)$ to the intervals $[0,+\infty)$ and  $(ab^{-1},0]$ respectively.
The following
 properties are elementary: 
\begin{enumerate}
\item[i)] $r^+ \in K_{a,b}^+, \ r^- \in K_{a,b}^-, \ Sr \equiv 0$ and $r(x,a,b) \to -a^2b^{-1}$ 
as $x\to +\infty$;
 \item[ii)] the inverse $\rho$ of $r$ is given by
$\rho(x,a,b)= ax/(a^2 +bx)$, and  $\rho''' (x)< 0$ for all $x>-a^2b^{-1}$;
\item[iii)] the equation $r(x,a,b)=-x$ has exactly two solutions: $x_1 =0$ and 
$x_{2}=(a+a^2)/b.$
\end{enumerate}
Futhermore, it can be proved that $r^+$ and $r^-$ are 
respectively the minimal element of $K_{a,b}^+$ 
and the maximal element of $K_{a,b}^-$ with 
respect to the usual order. The following 
slightly different result will play a key role in the sequel: 
\begin{lemma}
\label{cr} 
For all $y\in K_{a,b}$ with $(Sy)(x)<0$ for $x \not\in \{w: y'(w) = 0\}$,
we have   $r(x,a,b) < y(x)$ for all $x > 0$ 
and also $r(x,a,b) > y(x)$ for all $ x \in (ab^{-1},0)$.
\end{lemma}
\begin{proof} 
Take $g\in C^3(\R)$ and define $G(x) = g''(x)/g'(x)$ for all
$x \in D_g = \{x: g'(x)\neq 0\}$. Then we have $(Sg)(x) = G'(x) - (1/2)G^2(x), \ x \in D_g$.
Therefore, for every function
$g\in C^3(\R)$ with negative  Schwarzian, the associated 
function $G(x)$ satisfies the differential Riccati inequality 
$G'(x) - (1/2)G^2(x) < 0$  for all $x \in D_g$. 
Now, the lemma follows from standard comparison results 
(see, e.g. \cite[Theorem 5.III]{walter})
if we observe that $R=r''/r'$ and 
$Y=y''/y'$ satisfy
$(Sr)(x)=R'(x)-(1/2)R^{2}(x)=0$, \ $(Sy)(x)=Y'(x)-(1/2)Y^{2}(x) <
0$, for all $x\in (ab^{-1}, +\infty)\cap D_y$, and
$Y(0)=R(0)$.  Indeed, the above relations imply in cases (a)-(c) that
$R(x)>Y(x)$ for all $x>0.$  Now,  integrating $R$ and $Y$ over
the interval $(0,x)$, we get $r'(x) < y'(x)$.  Integrating $r'$
and $y'$ from $0$ to $x$  again, we obtain $r(x) < y(x)$ for
all $x>0$.
In cases (d)-(e) we obtain  using the above arguments  that $r(x) < y(x)$
for all $x\in (0,x^*)$. Now, since $r$ is strictly decreasing and $y$ reaches its minimum at 
$x^*$, it is obvious that the relation $r(x) < y(x)$ also holds for $x\geq x^*$.

Now the previous arguments allow  us to prove that  $r'(x) < y'(x)$  and 
$r(x) > y(x)$ for all $x \in (ab^{-1},0)$ in cases (a) and (d), where $y$ has no negative
inflexion points.

The proof  for (b) and (e) is slightly different if
the inflexion point $c_1$ of $y$ belongs to  the interval $ (ab^{-1},0)$. 
In this case, we can use the same arguments only to show that 
$r'(x) < y'(x)$ and $r(x) > y(x)$  for $x \in [c_1,0)$.
Next, by convexity arguments, 
$$
r(x)>r(c_{1})+r'(c_1)(x-c_1)> y(c_{1})+y'(c_1)(x-c_1)>y(x)
$$
for all 
$x \in (ab^{-1},c_1)$.

Finally,  case (c) can be studied analogously taking into account that $r$ is strictly
decreasing on $(ab^{-1},0)$, whereas $y$ reaches its maximum at $x^*$.
\end{proof}
As a by-product of the proof of Lemma \ref{cr},  we state the following corollary, which will be 
used in the proof of Theorem \ref{5678} when considering the case $f''(0)<0$.
\begin{corollary}
\label{bou}
Suppose that $f$ satisfies \textbf{(H)} and $f''(0) < 0$.
Then $f$ is bounded on $\R$. 
\end{corollary}
\begin{proof} 
Since $f$ satisfies \textbf{(H)}, the inequality 
$f''(0) < 0$ implies  that either $f''(x) < 0$ for all $x \leq 0$ (cases (f) and (g) in Fig. 1)
or $f$ has a global maximum at $x^*<0$ (see Fig. 1 (h)). Since in the latter case the
statement of the corollary is evident, we can assume that $f''(x)<0$ on $(-\infty ,0]$.

Next,  the function $g$ defined by $g(x)=-f(-x)$ satisfies $g'(0) = f'(0)$ and
$g''(x)=-f''(-x)>0$ for $x\geq 0$. Hence $g$ has not inflexion points on $[0, +\infty)$ and
we can use the property $(Sg)(x)= G'(x)-(1/2) G^2(x)<0$ for 
$x\in [0,+\infty)\cap D_g$ as it was done in the
first part of the proof of Lemma \ref{cr} to establish that $g(x) > 
r(x)=r(x,g'(0),g''(0)/2)$ for $x\in (0,+\infty)$. 
Since $r$ is strictly decreasing and $r(+\infty)=  2(f'(0))^2/f''(0) \in \R$, 
we can conclude that $g$ is bounded on $(0,+\infty)$. Thus $f(x)= -g(- x)$ 
is bounded on $(-\infty,0)$.
Finally, since $f$ satisfies \textbf{(H)}, $f$ is also bounded on $[0,+\infty)$.
\end{proof}
Now set $$a = f'(0), 2b = f''(0), \mu = ab^{-1} = 2f'(0)/f''(0),
r(x)=r(x,f'(0),f''(0)/2)$$ and define the continuous functions $A,B:
(\mu,+\infty)\to \R$ and $D:\R_+ \to \R$ by
$$
A(x)=x+r(x)+\frac{1}{r(x)}\int_{x}^{0}r(t)dt,
\quad 
B(x)=\frac{1}{r(x)}\int_{-r(x)}^{0}r(s)ds
\quad {\rm for \ } x\not=0, 
$$
$$
A(0) = B(0) = 0, \quad
D(x)= \left\{ \begin{array}{ll} A(x) & \textrm{if $r(x)<-x$},
\\ B(x) & \textrm{if $r(x)\geq -x$}. \end{array} \right.
$$
For the case $f'(0)<-1$, we will also use  the  function $R(x) = r(x, A'(0),
A''(0)/2)$ defined on the interval 
$(2A'(0)/A''(0), \infty) = (\nu,\infty)$.  Note that $A'(0)=f'(0)+1/2<0, 
\ A''(0) = f''(0)(1 + (6f'(0))^{-1})>0$. It is easy 
to check that $f'(0) < -1/6$ implies $(\nu,+\infty) \subset (\mu, +\infty)$, 
and that $A(x_2)=B(x_2)$, where $r(x_2) = -x_2 <0$. Also $B'(0)= -(r'(0))^2/2 = 
-(f'(0))^2/2$.  In the following lemmas we state other properties of the functions $A,B,D,R$.
\begin{lemma}
\label{31}
If $f'(0)< -1$ then
$A'(x)<0$ and $(SA)(x)<0$ for all $x \in (\mu,x_2)$.
\end{lemma}
\begin{proof} Since $f'(0) < -1$, we have  
$x_2 >0$ and $-xr(x) < r^2(x)$ for all $ x \in (\mu,x_2)\setminus\{0\}$. Hence  
$$\int_{x}^{0}r(t)dt \leq -xr(x) < r^{2}(x) \ {\rm and} \ 
A'(x)=r'(x)\bigg(1-\frac{\int_{x}^{0}r(t)dt}{r^{2}(x)}\bigg) < 0,  \ x \not=0.
$$
We get also $A'(0)=r'(0)+1/2<-1/2<0$.

Next, integrating by parts, we obtain
\begin{eqnarray*}
A(x)&=&r(x)+ \frac{xr(x)+\int_{x}^{0}r(t)dt}{r(x)}=r(x)+\frac{xr(x)+\int_{r(x)}^{0}vd\rho(v)}{r(x)}\\
   & =& r(x)+\frac{1}{r(x)}\int_{0}^{r(x)}\rho(v)dv=G(r(x)),
\end{eqnarray*}
where $\rho(v)=r^{-1}(v)$ and $ G(z)=z+\int_{0}^{1}\rho(vz)dv$. 

Then, by the formula for the Schwarzian derivative of the composition of two functions \cite{stri}, we obtain
 $$(SA)(x)=(SG)(r(x))(r'(x))^{2}+(Sr)(x)=(SG)(r(x))(r'(x))^{2}.$$
Hence, in order to prove that $SA<0$, we have to verify that $(SG)(r(x))<0$. On the other hand, 
$A'(x) <0$ if and only if $G'(r(x))>0$, so  
it suffices to show that $(SG)(r(x))<0$ when
$G'(r(x))>0$. Taking into account the fact that $\rho'''(z)< 0$ for $z=r(x)$, we have
$G'''(z)=\int_{0}^{1}v^{3}\rho'''(vz)dv <0$
and 
therefore
$
(SG)(r(x))<0.
$
\end{proof}
\begin{lemma}
\label{32}
If $f'(0) < -1$ then 
$(A(x)-R(x))x > 0$ 
for  $x \in (\nu,x_2) \setminus \{0\}$. 
\end{lemma}
\begin{proof} Recall that $R(x) = r(x,A'(0), A''(0)/2)$ and apply 
Lemmas \ref{cr}, \ref{31}.  
\end{proof}
Next, we will compare the functions $B(x)$ and $R(x)$ over 
the interval $[x_2,+\infty)$. In order to do that, we need 
the following simple result: 
\begin{lemma}
\label{new} If $(s,\zeta) \in \Pi = [-1,0] \times [-1.5, -1.25]$, then 
$$
L(\zeta,s) = s-\zeta -\ln(1+s -\zeta) +
\frac{2(\zeta+1/2)^2(s-\zeta)^2}{(2\zeta+1)\zeta^2+ (2/3)\zeta(s-\zeta)}< 0.
$$
\end{lemma}
\begin{proof} We have
$$
\frac{\partial L(\zeta,s)}{\partial s}= 
\frac{(s-\zeta)(s- A_-(\zeta))(s-A_+(\zeta))(12\zeta^2+16\zeta+3)}
{\zeta(1+s-\zeta)(6\zeta^2+\zeta+2s)^2},
$$
where 
\begin{eqnarray*}
A_{\pm}(\zeta) = -\frac{72\zeta^4+108\zeta^3+46\zeta^2+15\zeta +3}{2(12\zeta^2+16\zeta+3)}\\
\pm \frac{3(\zeta+1)(2\zeta+1)\sqrt{(2\zeta+1)
(72\zeta^3-12\zeta^2-6\zeta+1)}}{2(12\zeta^2+16\zeta+3)}.
\end{eqnarray*}
\noindent Note that $A_+$ and $A_-$ are continuous on the interval $J=[-\frac32,
-\frac23 -\frac{\sqrt{7}}{6})= [-1.5,-1.107..)$, and
$A_-(\zeta) < A_+(\zeta)$ for all
$\zeta
\in  J$. 
Now, it is straightforward to see (after several elementary 
transformations) that every root $\zeta$  of  $A_+(\zeta) = -1$ satisfies 
$\zeta(\zeta+1)^2(12\zeta^2 +16\zeta +3)(36\zeta^2 +12\zeta -7) =0$, 
and therefore belongs to the set 
$\{0,-1,-1/6 \pm \sqrt{2}/3, -2/3 \pm \sqrt{7}/6\}= \{0,-1,-0.638..., 
0.304..., -1.107...,-0.225...\}$. 
Since $A_+(-1.5) = -1.29991... < -1$, we obtain that 
$$ A_-(\zeta) < A_+(\zeta) < -1 \leq s \ {\rm for \ all \ } \zeta \in [-1.5,-1.107...)\,
,\, s\in [-1,0].$$ This implies immediately that 
$\partial L(\zeta,s)/\partial s < 0$ for all $\zeta \in [-1.5,-1.25]$ and $s\in
[-1,0].$  Thus $$\max_{(s,\zeta) \in \Pi}L(\zeta,s) = \max_{\zeta \in
[-1.5,-1.25]}L(\zeta,-1).$$  Finally, for $\zeta \in [-1.5,-1.25]$, we have 
$$\partial L(\zeta,-1)/\partial \zeta = 
-(36\zeta^2 +16\zeta -3)(\zeta+1)^3\zeta^{-2}(3\zeta+2)^{-2}(2\zeta-1)^{-2} > 0,$$
so that $\max_{\zeta \in [-1.5,-1.25]}L(\zeta,-1) = L(-1.25,-1) = -0.0006945... <0$. 
\end{proof}
\begin{lemma}
\label{45}
If $f'(0)  \in [-1.5, -1.25]$ then $B(x)>R(x)$
for all $x\geq x_2$. 
\end{lemma}
\begin{proof} Note that $B(x)= \tilde B(-r(x)),$ and
$R(x)= \tilde R(-r(x))$, where 
$$\tilde B(u)=\int_{0}^{1}r(zu)dz = 
\frac{\zeta}{u\theta}\big(u-\frac{1}{\theta}\ln(1+\theta u)\big), $$
and
$$ \tilde R(u)= R(\rho(-u))= - \frac{2(\zeta+1/2)^2u}{(2\zeta+1)\zeta+ (2/3)\theta
u},$$ with $\zeta=f'(0)$, $\theta = - f''(0)/(2f'(0)) >0$ and $\rho(x)=r^{-1}(x)$. 
Therefore  $B(x)-R(x) >0$ for $x \in [x_2, +\infty)$ if and only if 
$\tilde B(u) - \tilde R(u) > 0$ for all $u= -r(x) \in [x_2, -r(+\infty)) 
= \theta^{-1}[-\zeta-1, -\zeta)$. Finally, by Lemma \ref{new}, we have that 
$$\tilde B(u) - \tilde R(u) = \frac{\zeta}{u\theta^2}L(\zeta, \theta u + \zeta) > 0$$
for all $\zeta\in [-1.5,-1.25]$ and  $u \in \theta^{-1}[-\zeta-1, -\zeta)$, since 
in this case $s = \theta u +\zeta \in [-1, 0]$.
\end{proof}

Lemmas \ref{32} and  \ref{45} together yield the
\begin{corollary}
\label{cora}
If $f'(0)  \in [-1.5, -1.25]$, then $D(x) > R(x)$ for $x > 0$. 
\end{corollary}

\section{Proof of Theorem \ref{5678}}

In this section,  we prove Theorem \ref{5678}. Thus we assume that $f$ satisfies
all conditions {\bf (H)}.

Denote by $C$ the space of all continuous real functions $\gamma$ on the interval
$[-1,0]$, with the norm given by $\|\gamma\|=\max_{-1\leq t\leq 0}|\gamma (t)|$. Since
$f$ is continuous, each $\gamma\in C$ determines a unique solution
$x=x(\cdot,\gamma)$ to Eq. (\ref{wr}) on $[-1, +\infty)$ such that $x(t)=\gamma (t)$,
$t\in [-1,0]$.
 It is well-known that the application 
$F^{t}\gamma:\R_+\times C \to C$ given by $F^{t}\gamma(s)=x_t(s)=x(t+s,\gamma), 
\ s \in [-1,0],$
defines a continuous  semiflow on $C$. Moreover, this semiflow is point dissipative:
\begin{lemma}
\label{lem}
Assume that $f$ is bounded below and  satisfies \textbf{(H1)}.  Then there exists $K_0
>0$ such that for any $\gamma \in C$ we have 
$\limsup\limits_{t \to +\infty}\|F^t\gamma\| < K_0$. 
\end{lemma}
\begin{proof}
See for example \cite[Proposition 2.1]{mp}. 
\end{proof}
Now,  fix an arbitrary $\gamma\in C$.
It follows from Lemma \ref{lem} that 
the $\omega$-limit set $\omega(\gamma)$ of the 
trajectory $\{F^{t}\gamma: t \in \R_{+}\} \subset C$ 
is an  invariant and compact set.
Write $m=m_{\gamma} =\min_{\alpha \in \omega(\gamma)}\alpha(0)$
and 
$M= M_{\gamma} =\max_{\alpha \in \omega(\gamma)}\alpha(0)$. 

Evidently, Theorem \ref{5678} is proved if we demonstrate  
that $M_{\gamma} = m_{\gamma} =0$ once $-f'(0) \in [0,1.5]$. 
We only need to consider the case  $m < 0 < M$.
Indeed, otherwise (i.e. $m \geq 0$ or $M \leq
0$) every solution
$x(t,\alpha),\ \alpha\in \omega(\gamma),$ is bounded and monotone (due to ({\bf H1})), and the
limit values
$x_{\pm} = x(\pm\infty,\alpha) \in \R$  are steady states of (1.1).  Now, since
(1.1) possesses a unique equilibrium $x(t) \equiv 0$, we conclude that
either $0= x_- \leq x(t,\alpha) \leq  x_+ = 0$ or  $0= x_- \geq x(t,\alpha) \geq  x_+ = 0$ for
all $t\in\R$, so that $\omega(\gamma) = \{0\}$. Therefore,  $m = M =0$ if  $m \geq 0$ or $M \leq
0$, and in the sequel we will always assume that $m < 0 < M$.


In the following lemmas we establish some relations between $m$ and $M$ which are 
needed in the proof of Theorem \ref{5678}.
\begin{lemma}
\label{44}
We have $m > D(M)$ and $m > r(-r(M)/2)$.
\end{lemma}
\begin{proof} 
First we assume that $r(M) \leq -M$. Then 
$t_{1}=M/r(M) \in [-1,0]$. Next,  $x(t)=r(M)t,\ t \in [t_{1},t_{1}+1]$
is the solution of the initial value problem $x(s)=M, \ s \in [t_1-1,t_1]$
for
\begin{equation}
\label{ep}
x'(t)=r(x(t-1)). 
\end{equation}
Since $m = \beta(0) = \min_{\alpha \in \omega(\gamma)}\alpha(0)$ 
for some $\beta \in \omega(\gamma)$, and since $\omega(\gamma)$
is an invariant set, we obtain the existence of a solution $\tilde z(t):\R \to \R$ 
to (\ref{wr}) such that $\tilde z(s) = \beta(s), \ s \in [-1,0]$ and 
$\tilde z_t \in \omega(\gamma)$ for every $t \in \R$. 
Obviously, $z(t) = \tilde z(t-1)$ also satisfies (\ref{wr}) for all $t \in \R$ 
and is such that $z(t) \geq z(1)=m$ for all $t \in \R$. 
This implies $0 = z'(1)= f(z(0))$. Finally, by 
hypothesis  {\bf (H1)}, we get $z(0) =0$.  

Let fix now this solution $z = z(t)$. 
Clearly $M=x(t) \geq z(t)$ for all $t \in [t_{1}-1,t_{1}]$.
Moreover, we will prove that $x(t) \geq z(t)$ for all $t \in [t_{1},0]$.

Indeed, if this is not the case we can find 
$t_{*} \in [t_{1},0)$ such that $x(t_{*})=z(t_{*})$ and 
$x(t)\geq z(t)$ for all $t \in [t_1-1,t_{*}]$. We claim  that 
\begin{equation}
\label{zvx}
z'(t)> x'(t)\ \mbox{ for all }\  t \in [t_{*},0]. 
\end{equation}
We distinguish two cases:
if $z(t-1)>0$ then, using Lemma \ref{cr}, we obtain 
$
x'(t)=r(x(t-1))\leq r(z(t-1))<f(z(t-1))=z'(t).
$
Next if $z(t-1)\leq 0$, then $z'(t)=f(z(t-1))\geq 0>r(x(t-1))=x'(t)$.

 After integration 
 over $(t_*,0)$, and using $x(0)=z(0)=0$, it follows from (\ref{zvx}) that $z(t_*)<x(t_*)$,
which is a contradiction.

Thus $x(t)\geq z(t)$ for $t \in [t_{1},0)$ and, arguing as above, we obtain
 \begin{eqnarray*}
& &m = \int_{0}^{1}z'(s)ds=\int_{0}^{1}f(z(s-1))ds >\int_{0}^{1}r(x(s-1))ds=\\
& & \int_{-1}^{t_{1}}r(M)ds + \int_{t_{1}}^{0}r(x(u))du =
M + r(M)+  \int_{t_{1}}^{0}r(r(M)u)du = A(M).
\end{eqnarray*}
Now in the  general case (i.e. we do not assume that $r(M)\leq -M$), we will prove
that
\begin{equation}
\label{fzr}
f(z(t))>r(r(M)t)\, ,\,\forall\ t\in (-1,0).
\end{equation}
To do this, we first  show that $z'(t)>r(M)$ for all $t\in (-1,0).$ 
Indeed, if $z(t-1)>0$ then, by Lemma \ref{cr}, $z'(t)=f(z(t-1))>r(z(t-1)\geq r(M)$.
On the other hand, if $z(t-1)\leq 0$ then $z'(t)=f(z(t-1))\geq 0>r(M)$. Hence, 
$$z(t)=-\int_t^0z'(s)ds<-\int_t^0r(M)ds=r(M)t\, ,\, t\in (-1,0).$$
To obtain (\ref{fzr}) we only have to note that $f(z(t))>r(z(t))>r(r(M)t))$ if
$z(t)>0$ and $f(z(t))\geq 0>r(r(M)t))$ if $z(t)\leq 0$.

Now, using (\ref{fzr}), we obtain
$$
m=z(1) =\int_{0}^{1}f(z(s-1))ds > \int_{0}^{1}r(r(M)(s-1))ds  = B(M) .
$$ 
Finally, applying Jensen's integral inequality (see \cite[p.~110]{roy}) , we have  
$$m>B(M) =\frac{1}{r(M)}\int_{-r(M)}^{0}r(s)ds \geq  r(-r(M)/2). $$
This completes the proof.
\end{proof}
As a consequence of Lemma \ref{45} and  Lemma \ref{44},
we obtain that $R(m),r(m)$ and $r(r(-r(M)/2))$ are 
well-defined and that $R(\nu,+\infty) \subset (\nu,+\infty)$ for  suitable values
of $f'(0)$:
\begin{corollary}
\label{cy}
\hspace{-2mm}We have $m > \mu, r(r(-r(M)/2)) > \mu$ if 
$f'(0) \in [-1.5,0)$ and $m > \nu$ if $f'(0)\in [-1.5,-1.25]$. 
\end{corollary}
\begin{proof} Indeed, for  $f'(0) \in (-2,0)$ and $f''(0) > 0$, we have 
$$
m > r(-r(M)/2) > r(-\frac{r(+\infty)}{2})= \frac{(f'(0))^3}{f''(0)(1-f'(0))} \geq 
\frac{2f'(0)}{f''(0)} = \mu.
$$ 
Next, $-A'(0) = -(f'(0)+0.5) \leq 1$ for $f'(0) \geq -1.5$, and 
Lemmas \ref{45}, \ref{44} lead to the estimate
$$m  > D(+\infty) = B(+\infty) \geq R(+\infty) = -A'(0)\nu \geq \nu.$$
This proves the corollary.
\end{proof}
\begin{lemma}
\label{43} Let $f'(0) \in [-1.5,0)$. We have
$M < r(m)$.  Moreover, if $f'(0)\in [-1.5,-1.25]$ then
$M < R(m)$. 
\end{lemma}
\begin{proof} We have that $r(m)$ is well defined and $[m, +\infty) \subset [\mu, +\infty)$
since $f'(0) \in [-1.5,0)$ (see Corollary \ref{cy}). 
Take now $\theta \in \omega(\gamma)$ such that $y(t)=y(t,\theta)$ satisfies $y(1)=M$ 
and, consequently, $y'(1)= y(0)=0$. First we prove that $f(y(s))<r(m)$ for  $s\in [-1,0]$ and $f'(0)\in [-1.5,0)$.
Indeed, we have  {$f(y(s))<r(y(s))\leq r(m)$} if $y(s)<0$ and $f(y(s))\leq 0<r(m)$ if
$y(s)\geq 0$. Thus
$$
M=y(1)=\int_0^1f(y(s-1))ds<\int_0^1r(m)ds=r(m).
$$
Now, if $f'(0)<-1$ then  $r(m)>-m$, from which it follows that $t_2=m(r(m))^{-1} \in (-1,0]$. Next,
$x(t)= r(m)t$, with $\ t \in [t_2,t_2+1]$, is the solution of the initial value 
problem $x(s)=m,\ s \in [t_2-1, t_2]$ for Eq. (\ref{ep}). 
Now we only have to argue as in the proof of Lemma \ref{44}
to obtain the inequality $M<A(m)$.
Finally, by Lemma \ref{32} and Corollary \ref{cy}, we obtain  
$M<A(m)<R(m)$ when $f'(0)\in [-1.5,-1.25]$. 
\end{proof}
\begin{remark} 
\label{remix}
Assume that $(f(x) - f'(0)x)x >0$ 
for all $x \neq 0$.  Replacing $r(x)$ with $r_1(x)= f'(0)x$ in 
the proof of Lemma \ref{44}, we can observe that it still 
works without any change if we set $A(M) = (f'(0)+1/2)M$ and 
$B(M) = - (f'(0))^2M/2$. The same observation is valid 
for the proof of Lemma \ref{43} (up to the last 
sentence beginning from the word ``Finally"). 
Therefore, under the above assumption, 
we have $m > A(M) =  (f'(0)+1/2)M$ and $M < A(m) =(f'(0)+1/2)m$ 
once $f'(0) \leq -1$. Also
$m > B(M)= - (f'(0))^2M/2$ and $M < r_1(m) = f'(0)m$ if $f'(0) < 0$. 
\end{remark}
\begin{proof} [{\it Proof of Theorem \ref{5678}}]  We will reach  
a contradiction if we assume that $m<0<M$.
Suppose first that $f''(0)>0$.
If $f'(0) \in (-1.5,0)$, in view of Lemmas \ref{44}, \ref{43} and Corollary \ref{cy}
we obtain that $M < r(m) \leq r\circ r(-r(M)/2) = \lambda (M)$ with the rational 
function $y =\lambda(x)$. Now, $\lambda(M) < M$ for $M > 0$ 
if $\lambda'(0) =(1/2) |f'(0)|^3< 1$. Therefore,  if $f'(0) \in [-1.25, 0)$ we obtain 
the desired contradiction under the assumption $M>0$. 

Now let $f'(0)\in [-1.5, -1.25]$ and, consequently, 
$R'(0) = f'(0) + 0.5 \in [-1,-0.75]$.
In this case Corollary \ref{cora},  
Lemma \ref{44} and Lemma \ref{43} imply that $M < R(R(M))$.
As $R\circ R(x) \leq  x$ for all $x > 0$ whenever 
$(R\circ R)'(0) = (R'(0))^2 \leq 1$, we obtain a contradiction again.
Therefore the solution $x(t)\equiv 0$ of Eq. (\ref{wr}) 
is globally attracting if $f''(0) > 0$ and  $f'(0) \in [-1.5,0)$. 

Assume now that $f''(0) < 0$. The change of variables $y(t)= -x(t)$ 
transforms (\ref{wr}) into $y'(t) = g(y(t-1))$ with $g(x)= -f(-x)$. 
It is easily seen that $g''(0) > 0$ and that
$g$ satisfies all properties from \textbf{(H)} (note that by 
Corollary \ref{bou}, $g$ is  bounded below).
Applying 
now the first part of the proof to the modified equation $y'(t) = g(y(t-1))$,
we reach the same contradictions if we assume that $m<0<M$. 
Hence our theorem is also proved when $f''(0) < 0$.
Note that this change of variables transforms the cases (f), (g) and (h) from
Fig. 1 into (b), (c) and (e) of the same figure respectively.

Finally, take $f''(0) = 0$. In this case, $x =0$ is an
inflexion point for $f$ and  $(f(x)-f'(0)x)x > 0$ if 
$x \not=0$. Therefore it is natural to employ the linear 
function $r_1(x)=f'(0)x$  instead of the rational function 
$r(x)$ used in the case $f''(0)>0$.
Now, by Remark \ref{remix} we have that 
$m > A(M)=(f'(0)+1/2)M$ and 
$M < A(m)=(f'(0)+1/2)m$ if $f'(0)\in [-1.5,-1)$. Hence
$M < (f'(0)+1/2)^2M \leq M$, a contradiction. Let now $f'(0)\in [-1,0)$.
By the same Remark \ref{remix} we obtain that
$m > B(M)=-(1/2) (f'(0))^2 M > -(1/2) (f'(0))^3 m$.
Thus $-(f'(0))^3/2 >1$, a contradiction.
\end{proof}

\section*{Acknowledgements}

We thank the referee of  this paper for his/her careful and
insightful critique.

\end{document}